\documentstyle[a4,12pt]{article}
\begin{document}
\begin{center}
{\bf ON THE PRODUCTS IN THE FINITE GROUPS }\\
\vspace{0.5cm}
{\bf V.V.Genk}\\
\vspace{0.5cm}
{\small\it Belarusian State Politechnical College \\
Tsentralnaya sq.2, Molodechno, Minsk region, 222310, The Republic of Belarus \\
 E-mail: {\rm gravitation@mail.admiral.ru}}
 \end{center}
\vspace{0.5cm}
{\small\begin{center}
\begin{tabular}{p{13cm}}
     \quad All possible products of all elements of an odd order finite group are
     considered. A set  of all such products is called as a $K$-set.

     \quad A hypothesis of $K$-set coincidence of any group of an odd order with its
     commutant is proposed and the hypothesis validity for groups with a commutant
     of a simple structure is shown.
\end{tabular}\end{center} }

\vspace{0.5cm}

   Sequences of elements of a finite group can
form ``a good" in this or that sense, part of
a group (see, for example, \cite{2}-\cite{7}), and a set of
products of all commutators is a commutant of
a group. A question arises: what information
about a group does a set made of every possible
products of all elements of the finite group
taken in an unspecified order carry?

\medskip
   Let\,  $G = \{1, a_1, {a_1}^{-1},\ldots, a_n, {a_n}^{-1}\}$\, be a finite
group of an odd order\, $2n+1$ (later on simply
{\it``an odd group"}), where\, $a_{i}^{\pm1}$ --- the whole of
different elements of the group.

   Let's form a set:
   \[K=\{g_{\sigma (1)}\ldots g_{\sigma(2n)}\mid \sigma\in S_{2n},\quad g_i\in G\setminus1,
   \quad g_i\neq g_j\quad \mbox{for } i\neq j \},\]
where \, $S_{2n}$ --- a set of substitutions in\, $2n$\, symbols.
Let's later on call a set $K$ as a $K$-{\it set}.
This paper proposes a hypothesis on the structure
of a $K$-set of odd groups and a special case of
this hypothesis is proved.

\medskip
   \underline{{\it T h e\quad b a s i c\quad h y p o t h e s i s}.}

   \quad{\bf A  $K$-set  of  any  odd  group  coincides  with  its  commutant}.

\medskip
  The following two suggestions take place
(the first is evident, the second is proved by
induction on word length out of a $K$-set).

\medskip
   \underline{{\it S u g g e s t i o n\quad1}.}

   \quad{\sl   A  $K$-set  of  an  odd  group  is  its
invariant  subset,  containing  the  unit  of  a  group  and  also  together
with  any  element  and  its  inverse}.

\medskip
   \underline{{\it S u g g e s t i o n\quad2}.}

 \quad{\sl A  K-set  of  an  odd  group  is contained  in  its  commutant}.

\medskip
  Just from these suggestions by the Feit-Thompson
theorem such a  consequence  arises (indicated by prof. L.S.Kazarin):

\medskip
    \underline{{\it C o n s e q u e n c e}.}

    \quad{\sl A  K-set  of  an  odd  group  is
solvable  too,  i.e.  a  repeat  procedure  of  taking  of  a K-sets
through  the  finite  number  of  steps  leads  to  the  unit  of  a  group}.

\medskip
  When proving the following theorem commutator
identities are systematically used
\begin{equation}\label{1}
  [a,b]=[a,ab],\quad [a,b]=[b^{-1},ab],\quad [a,b]^{-1}=[b,a],
\end{equation}
and an obvious  lemma  as well:

\medskip
    \underline{\it{L e m m a}.}

    \quad{\sl Not  any  other element  of  an  odd  group  but
    the  unit  can  be  conjugated  with its  inverse.}

\medskip
    \underline{\it{T h e o r e m}.}

    \quad{\sl The  product  of  two  commutators  in  an  odd
group  belongs  to  its  $K$-set.}

\medskip
  \noindent {\it P r o o f.}

     It's sufficient to check that the
inclusion takes place when one of the commutators
is made of elements which themselves or their
inverse enter the other commutator. Sixteen
variants of accurance of common elements or their
inverse into commutators are possible with an
accuracy to designation of the elements of a group:
\begin{equation}\label{2}
\begin{array}{cccc}
    {[a_1,a_2][a_1,a_3],} & {[a_1,a_2][a_2,a_3],} &{[a_1,a_2][a_3,a_2],} \\
    {[a_1,a_2][a_3,a_1],} & {[a_1,a_2][a_3,a_1^{-1}],} &{[a_1,a_2][a_1^{-1},a_3],}\\
    {[a_1,a_2][a_2^{-1},a_3],}&{[a_1,a_2][a_3,a_2^{-1}],}&{[a_1,a_2][a_2,a_1],}\\
    {[a_1,a_2][a_1,a_2],}&{[a_1,a_2][a_1,a_2^{-1}],}&{[a_1,a_2][a_1^{-1},a_2],}\\
    {[a_1,a_2][a_2,a_1^{-1}],}&{[a_1,a_2][a_2^{-1},a_1],}&{[a_1,a_2][a_2^{-1},a_1^{-1}],}\\
    {[a_1,a_2][a_1^{-1},a_2^{-1}].}
 \end{array}
\end{equation}

\medskip
  Let's consider the first of them in more detail.
The rest are analyzed similarly or reduced to
already analyzed ones. It's clear that the product
of elements\, $a_1a_2\notin~\{a_1,a_2,a_1^{-1},a_2^{-1}\}$.

   Let
    \[a_1a_2\in~\{a_3,a_3^{-1},a_i^{\varepsilon}\},
    \quad \mbox{where } i>3,\quad \varepsilon=\pm1.\]
  1) If\,\, $a_1a_2=a_i^{\varepsilon},$\\
  then by the second of the identities (1):
  \[ [a_1,a_2][a_1,a_3] = [a_2^{-1},a_i^{\varepsilon}][a_1,a_3]\in K.\]
  2) Let\,\, $a_1a_2=a_3$,\,\,then:
  \[[a_1,a_2][a_1,a_3] = [a_2^{-1},a_3][a_1,a_3].\]
  Consider the product of elements\,\, $a_2^{-1}a_3$. As shown above
      \[a_2^{-1}a_3\in\{a_1,a_1^{-1},a_i^{\varepsilon}\},\quad i>3.\]
 {\it a})\,\, For\, $a_2^{-1}a_3=a_i^{\varepsilon}$\, by the
    first of the identities (1):
      \[[a_2^{-1},a_i^{\varepsilon}][a_1,a_3]\in K.\]
 {\it b})\,\, If\, $a_2^{-1}a_3=a_1$,\, then\, $a_3a_1=a_1a_2a_2^{-1}a_3=a_1a_3$\, and
\[[a_1,a_2][a_1,a_3]\in K.\]
 {\it c})\,\, It remains\, $a_2^{-1}a_3=a_1^{-1}$,\, from which\,
$a_1^{-1}=a_2^{-1}a_1a_2$,\, i.e. $a_1\sim a_1^{-1}$\, and this case is
impossible due to the lemma.

\noindent 3) Let, finally,\,\, $a_1a_2=a_3^{-1}$. Then\\
\[[a_1,a_2][a_1,a_3] = [a_2^{-1},a_3^{-1}][a_1,a_3].\]
 Considering the product
 \[a_2^{-1}a_3^{-1}\in\{a_1,a_1^{-1},a_i^{\varepsilon}\},\quad i>3 \]
as above one can make sure that in any case
\[[a_1,a_2][a_1,a_3]\in K.\]
   All another variants in (2) may be considered similarly. So the
    proof is completed.

\medskip
   Thus, the basic hypothesis is true for odd groups
with a commutant elements of which reach the limit
of products of not more than two commutators.
What concerns noncommutative groups of an even
order the question about a structure of $K$-sets
of such groups remains open. An example of
symmetric group $S_3$ shows that a $K$-set of an even
group can disagree with its commutant.

\medskip
   The author expresses his thanks to prof. L.S.Kazarin
for useful discussions and valuable remarks.

\vspace{0.5cm}
\end{document}